\documentstyle[12pt]{article}
\def\lqq{\lq \lq }

\begin{document}
 
\begin{titlepage}

\rightline{Preprint  ICN-UNAM/96-02}
\rightline{January 1996}
\begin{center}
{\large Interesting relations in Fock space$^{\star}$
\vskip 0.5cm
by 
\vskip 0.5cm
 Alexander Turbiner}$^{\dagger}$
\vskip 0.5cm
{\it Instituto de Ciencias Nucleares, UNAM, A.P. 70-543, 04510 M\'exico D.F}
\end{center}

\vskip 0.5 cm
\begin{center}
{\large ABSTRACT}
\end{center}
\vskip 0.5 cm
\begin{quote}
Certain non-linear relations between the generators of the ($q$-deformed)
Heisenberg algebra are found. Some of these relations are invariant 
under quantization and $q$-deformation.
\end{quote}
\vfill
\noindent
$^{\star}$ Invited talk at International Conference 
on Complex and Hypercomplex Analysis, Mexico-city, Mexico (1994);
to be published in Proceedings

\noindent
$^\dagger$ On leave of absence from the Institute for Theoretical
and Experimental Physics, Moscow 117259, Russia \\
E-mail: turbiner@nuclecu.unam.mx or turbiner@axcrna.cern.ch
\end{titlepage}

In this Talk I want to present different relations appearing
in the Fock space generated by the q-deformed Heisenberg algebra.
In a certain particular case, some of these relations can be 
summarized by the theorem:

{\bf \it  (main) THEOREM}. The following differential identities hold 
for any natural number $n$ and $\delta \in {\bf C}$
\begin{equation}
\left(x^2{d \over dx}-n x\right)^{n+1} \ =\  x^{2n+2} {d^{n+1}\over dx^{n+1}}\ ,\quad 
\left({d^2 \over dx^2} x-n {d \over dx} \right)^{n+1} \ =\  {d^{2n+2} \over dx^{2n+2}} x^{n+1}\ ,
\end{equation}
\begin{equation}
\left(x{d \over dx} x\right)^n \ =\  x^n {d^n \over dx^n} x^n\ ,\quad 
\left({d \over dx} x {d \over dx}\right)^n \ =\ {d^n \over dx^n} x^n {d^n \over dx^n}\ ,
\end{equation}
\begin{equation}
x{d \over dx}\left(x{d \over dx}-1\right) \ldots\left(x{d \over dx}-n\right)=x^{n+1}{d^{n+1} \over dx^{n+1}}\ ,
\end{equation}
\begin{equation}
\prod_{k=0}^n
\left[x\left(1-e^{-\delta {d \over dx}}\right)-k\delta\right] \ =\  x^{(n+1)} 
\left(1-e^{-\delta {d \over dx}}\right)^{n+1}\ ,
\end{equation}
where $x^{(n+1)}\ =\ x(x - \delta) (x - 2\delta) \ldots (x - n\delta)$.

The proof can be carried out by induction.
Later we will show that these relations can be easily generalized to 
various Fock spaces. 
 
Take some 3-dimensional complex algebra with an identity operator and with 
two elements $a$ and $b$, obeying the relation 
\begin{equation}
\label{e1}
ab - q ba \ =\ p\ , 
\end{equation}
where $p, q$ are any complex numbers. The algebra with the identity operator 
generated by the elements $a,b$ and obeying (1) is usually called the $q$-deformed  
Heisenberg algebra $h_q$ (See, for example, papers \cite{1} and references 
and a discussion therein).  
The parameter $q$ is called the parameter of quantum deformation. If $p=0$,
then (5) describes the non-commutative (quantum) plane (See, for example, \cite{2}), while
if $q=1$, then $h_q$ becomes the ordinary (classical) Heisenberg algebra 
and the parameter 
$p$ plays a role of the Planck constant (See, for example, \cite{3}). The operator linear 
space of all holomorphic functions in $a,b$ with vacuum
\[
a|0>\ =\ 0 
\]
is called a {\it Fock} space.
Then  the following theorem holds:

{\bf \it THEOREM 1 \cite{4}.}
For any $p,q \in {\bf C}$ in (\ref{e1}) the following identities are true
\renewcommand{\theequation}{6.{\arabic{equation}}}
\setcounter{equation}{0}
\begin{equation}
\label{e2.1}
(aba)^n\ = \ a^n b^n a^n , \quad n=1,2,3,\ldots \ ,
\end{equation}
and, if $q \neq0$, also
\begin{equation}
\label{e2.2}
(bab)^n\ = \ b^n a^n b^n , \quad n=1,2,3,\ldots \ .
\end{equation}
\renewcommand{\theequation}{\arabic{equation}}
\setcounter{equation}{6}
The proof can be carried out  by induction using the following easy lemma

{\bf \it  LEMMA 1}. 
For any $p,q \in {\bf C}$ in (\ref{e1}),  
\renewcommand{\theequation}{7.{\arabic{equation}}}
\setcounter{equation}{0}
\begin{equation} 
\label{e3.1}
a b^n \  -\ q^n b^n a \ = \ p \{ n \} b^{n-1} \ ,
\end{equation}
\begin{equation} 
\label{e3.2}
a^n b \  -\ q^n b a^n \ =\  p \{ n \} a^{n-1} \ ,
\end{equation}
where $n=1,2,3 \ldots $ and $\{ n \}= {1-q^n \over 1-q}$ is the so-called $q$-number
(See, for example, \cite{5}).

\renewcommand{\theequation}{\arabic{equation}}
\setcounter{equation}{7}

Now let us proceed to the proof of (6.1). For $n=1$ relation (6.1) holds trivially.
Assume that (6.1) holds for some $n$ and check that it holds for $n+1$. 
It is easy to write the chain of equalities :
\[
(aba)^{n+1}\ =\ (aba)^n aba \ =\ a^n b^n a^n aba\ 
=\ a^n (b^na)(a^nb)a\ =
\]
\[
=\ a^n ({ab^n \over q^n} - {p \over q^n} \{ n \} b^{n-1}) 
(q^n ba^n + p \{ n \} a^{n-1})a \ =\ a^{n+1} b^{n+1} a^{n+1} \ ,
\qquad (\diamondsuit)
\]
where the relations (7.1) and (7.2) were used in this chain. In an analogous way, one can prove (6.2). {\it q.e.d.}

 Theorem 1 leads to two corollaries, both of them easily verified:

{\bf \it Corollary 1}. 
For any $p,q \in {\bf C}$ in (\ref{e1})  and natural numbers 
$n,k$, 
\begin{equation}
\label{e4}
(\underbrace{ ababa....aba}_{2k+1})^n = 
\underbrace{a^n b^n a^n ...b^n a^n}_{2k+1}  , \quad n,k=1,2,3,\ldots
\end{equation}

{\bf \it Corollary 2}. 
Let  $T^{(n)}_k= \underbrace{a^n b^n a^n \ldots b^n a^n}_{2k+1} $.
Then, for any $p,q \in {\bf C}$ in (\ref{e1}),  the following relation holds
\renewcommand{\theequation}{9.{\arabic{equation}}}
\setcounter{equation}{0}
\begin{equation}
\label{e5.1}  
 [ T^{(n)}_k , T^{(m)}_k ]\ = \ 0 \ ,
\end{equation}
as well as the more general relation:
\begin{equation}
\label{e5.2} 
  [ T^{(n_1)}_kT^{(n_2)}_k \ldots T^{(n_i)}_k ,  
   T^{(m_1)}_k T^{(m_2)}_k\ldots T^{(m_j)}_k]\ = \ 0 \ ,
\end{equation}
\noindent
where $[\alpha, \beta] \equiv \alpha\beta - \beta \alpha$ is the standard commutator 
and  $\langle n \rangle,\langle m \rangle$ are sets of any non-negative, 
integer numbers.

\renewcommand{\theequation}{\arabic{equation}}
\setcounter{equation}{9}

It is evident that formulas (8), (9.1) and (9.2) remain correct under  
the replacement
$a \rightleftharpoons  b$, if  $q\neq 0$.  It is worth noting that 
the algebra $h_q$ under an appropriate choice of the parameters $p,q$ 
has a natural representation
\[
a = x\ ,\ b = D  \qquad  (\star)
\]
where the operator  $D f(x)={f(x)-f(qx) \over x(1-q)}$ is the dilatationally-invariant shift 
operator and it is usually called the Jackson symbol 
(See, for example, \cite{5}). The relation (6.1) then becomes
\[
(xDx)^n \ =\ x^n D^n x^n \ ,
\]
while relation (6.2) becomes
\[
(DxD)^n \ =\ D^n x^n D^n \ .
\]

Since for $q \rightarrow1$,  the operator $D \rightarrow {d\over dx}$, 
the above relations become differential identities (2).

{\bf \it THEOREM 2 \cite{6}.}
For any $p,q \in {\bf C}$ in (\ref{e1})  and natural $n,m$  the following identities  hold:
\renewcommand{\theequation}{10.{\arabic{equation}}}
\setcounter{equation}{0}
\begin{equation}
\label{e6.1}
[a^nb^n, a^mb^m]\ =\ 0 \ ,
\end{equation}
\begin{equation}
\label{e6.2}
[a^nb^n, b^ma^m]\ =\ 0 \ ,
\end{equation}
\begin{equation}
\label{e6.3}
[b^na^n, b^ma^m]\ =\ 0 \ .
\end{equation}
\renewcommand{\theequation}{\arabic{equation}}
\setcounter{equation}{10}

The proof  is based on the following easy lemma

{\bf \it LEMMA 2.}
For any $p,q \in {\bf C}$ in (\ref{e1})  and natural number $n$  
\renewcommand{\theequation}{11.{\arabic{equation}}}
\setcounter{equation}{0}
\begin{equation}
\label{e7.1}
a^n b^n \ =\ P(ab) \ ,
\end{equation}
\begin{equation}
\label{e7.2}
b^n a^n \ =\ Q(ab) \ ,
\end{equation}
where $P,Q$ are some polynomials in one variable of order not higher than $n$.

\renewcommand{\theequation}{\arabic{equation}}
\setcounter{equation}{11}

Let us introduce a notation $t^{(n)}=a^n b^n $ or  $b^n a^n$, in which the order of the 
multipliers is not essential. The statement of the Theorem 2 can now be written 
as $[t^{(n)},t^{(m)}]=0$ and the following is true;

{\bf \it Corollary 3}. 
The commutator
\begin{equation}
\label{e8} 
 [ t^{(n_1)} t^{(n_2)} \ldots t^{(n_i)}, \  t^{(m_1)} t^{(m_2)} \ldots t^{(m_j)}]\ = \ 0 
\end{equation}
holds for any  $p,q \in {\bf C}$ in (\ref{e1}) and any sets $\langle n \rangle, \langle m \rangle$  
of non-negative, integer numbers.

One can make sense of (6.1) and (6.2), (8), (9.1) and (9.2), (10.1)-(10.3), 
(12) as follows: In the algebra of 
polynomials in $a,b$ there exist relations invariant under a variation of the parameters 
$p,q$ in (5). Also formulas (6.1) and (6.2), (8) can be interpreted as formulas of a certain
special ordering other than the standard lexicographical one. 

A natural question can be raised: Is the existence of the relations (6.1) 
and (6.2), (8), (9.1) and (9.2), 
(10.1)-(10.3), (12) connected unambiguously to the algebra $h_q$, 
or are there more general algebra(s) leading to those relations? 
Some answer is given by the following theorem

{\bf \it THEOREM 3 \cite{6}.}
If two elements $a,b$ of a Banach algebra with unit element are related by
\begin{equation}
\label{e9}
ba = f(ab) \ ,
\end{equation}
where $f$ is a holomorphic function in a neighbourhood of 
$\hbox{Spec}\,\{ab\} \cup \hbox{Spec}\,\{ba\} $, then relations (6.1), 
(8), (10.3) hold. If, in addition, the function $f$ is single-sheeted, 
then (10.2) also holds.

The proof is essentially based on the fact that, if the function $f$ 
in (13) is holomorphic, then for {\it any} holomorphic $F$
\begin{equation}
\label{e9a}
b F(ab)\ =\ F(ba) b \ ,
\end{equation}
which guarantees the correctness of the statement (11.2) of Lemma 2, 
although $Q$ is no longer polynomial. This immediately proves (10.3). 
The relation (6.1) can be proved by induction and an analogue of the 
logical chain ($\diamondsuit$) is
\[
(aba)^{n+1}\ =\ (aba)^n aba \ =\ a^n b^n a^n aba \ =\ 
a^n (b^na^n) aba\ =\ 
\]
\[
= \ [ b^na^n = Q (ab), {\hbox{see Lemma 2}}]\ = 
\]
\[
 =\ a^n Q(ab) (ab) a \ =\ a^n (ab) Q(ab) a \ =\  
a^{n+1} b^{n+1} a^{n+1} \ .
\]
An extra condition that $f$ be single-sheeted implies that 
$ab = f^{-1}(ba)$, which immediately leads to the statement (10.2). 
{\it q.e.d.}

It is evident that the replacement $a \rightleftharpoons b$ in 
Theorem 3 leads to the fulfilment of the equalities (6.2), (10.1) 
and (10.2) as well.

There exists another type of relations in Fock space stemming from 
the fact that some algebras are contained in the Fock space and 
these algebras can possess finite-dimensional representations. 
One can prove the theorem:

{\bf \it THEOREM 4}.
For any $p,q \in {\bf C}$ in (\ref{e1}) and $n=1,2,3,\ldots$ 
the following identities are true
\renewcommand{\theequation}{15.{\arabic{equation}}}
\setcounter{equation}{0}
\begin{equation}
\label{e14.1}
(b^2a-\{ n\} b)^{n+1}=q^{n(n+1)}b^{2n+2}a^{n+1}\ , 
\end{equation}
and, if $q \neq 0$, also
\begin{equation}
\label{e14.2}
(ba^2-\{ n\} a)^{n+1}=q^{n(n+1)} b^{n+1} a^{2n+2}\ , 
\end{equation}
where $\{ n \}= {1-q^n \over 1-q}$ is the $q$-number.

\renewcommand{\theequation}{\arabic{equation}}
\setcounter{equation}{15}

In order to prove this theorem, we need, at first, to state 
the following observation 

{\bf \it  LEMMA 3 \cite{ot}.} 
For any $p,q \in {\bf C}$ in (\ref{e1}) and $\alpha \in {\bf C}$, 
three elements of the Fock space 
\[
j^+\ =\ b^2a- \{\alpha\} b \ ,
\]
\begin{equation}
\label{e15}
j^0\ =\ ba- {\{\alpha\}\{\alpha+1\}\over \{2\alpha+2\}} \ ,
\end{equation}
\[
j^-\ =\ a \ ,
\]
(modified by some multiplicative factors) obey $q$-deformed commutation 
relations
\[ 
q \tilde  j^0\tilde  j^- \ - \ \tilde  j^-\tilde  j^0 \
= \ - \tilde  j^-  \ ,
\]
\begin{equation}
\label{e16}
 q^2 \tilde  j^+\tilde  j^- \ - \ \tilde  j^-\tilde  j^+ \
= \ - (q+1) \tilde  j^0 \ ,
\end{equation}
\[ 
\tilde  j^0\tilde  j^+ \ - \ q\tilde  j^+\tilde  j^0 \ = \  \tilde  j^+  \ ,
\]
forming the algebra $s\ell_{2q}$. If $q \rightarrow 1$, 
these commutation relations become the standard ones for $s\ell_2$. 
If $\alpha=n$ is a non-negative integer number, the generators (16) 
form the finite-dimensional 
representation corresponding to the (operator) finite-dimensional 
representation space
\begin{equation}
\label{e17}
V_n=<1,b,b^2,\ldots, b^n>
\end{equation}
in the Fock space.

Validity of this Lemma can be checked by direct calculation. The proof 
of Theorem 4 is based on an evident fact \cite{7} that a positive-root (negative-root) generator in finite-dimensional representation taken in the 
power of the dimension of the finite-dimensional representation 
annihilates the space of the finite-dimensional representation and, correspondingly, acts in its complement. The operator $(j^+)^{n+1}$ at $\alpha=n$ annihilates (18) and hence it must be proportional 
(from the right) to $a^{n+1}$. Since $j^+$ is graded with the grading 
equals to (-1), $(j^+)^{n+1}$ must be proportional (from the left) to $b^{2n+2}$. What is left in this consideration is a value of possible multiplicative factor appearing in r.h.s. (15.1).
This factor is equal to $q^{n(n+1)}$ and can be found by direct 
calculation. Equation (15.2) can be proved analogously, with only minor modifications.

{\bf \it THEOREM 5 \cite{4}.}
For any $p,q \in {\bf C}$ in (\ref{e1}) and $n=1,2,3,\ldots$ the following identity holds :
\begin{equation}
\label{e19}
ba(ba-\{1\})\ldots (ba-\{ n\})\ =\ q^{n(n+1)\over2} b^{n+1}a^{n+1}
\end{equation}
The proof can be carried out by induction.
Now let us consider how the identity (\ref{e19}) looks for different representations of the algebra (5).

Taking the representation $(\star)$ for the algebra (5) and plugging 
it into (19), we arrive at an identity for the dilatationally-invariant 
shift operator $D$ : 
\begin{equation}
\label{e20}
xD(xD-\{1\})\ldots (xD-\{ n\})\ =\ q^{n(n+1)\over2} x^{n+1}D^{n+1}
\end{equation}
If $q \rightarrow 1$ the representation $(\star)$ degenerates into
\[
a= {d \over dx}\ ,\ b=x \qquad (\star \star)
\]
and (19) coincides to the identity (3). Recently a more general
representation of the algebra (5) was found then $(\star \star)$ 
at $q=1$ is characterized 
by a free parameter $\delta \in C$ \cite{st}:
\renewcommand{\theequation}{21.{\arabic{equation}}}
\setcounter{equation}{0}
\begin{equation}
\label{e22.1}
a={ (e^{ \delta {d \over dx}} - 1) \over \delta}\ ,
\ b=x e^{-\delta {d \over dx}} \ ,
\end{equation}
or, 
\begin{equation}
\label{e22.2}
a\ =\ {\cal D}_+ \ ,\ b\ =\ x(1 - \delta {\cal D}_-)\ ,
\end{equation}
where ${\cal D}_{\pm} f(x) \ =\ {f(x\pm\delta) - f(x) \over \pm\delta}$ 
are the (translationally-invariant) finite-difference operators.
After substitution of (21.1) in (19) for $q=1$ a simple transformation  
the differential identity (4) appears. This identity takes a slightly different form
 if the realization (21.2) is used :
\renewcommand{\theequation}{\arabic{equation}}
\setcounter{equation}{21}
\begin{equation}
\label{e23}
xD_-(xD_--1)\ldots (xD_--n)\ =\ {\delta}^{n+1} x^{(n+1)}D_-^{n+1}
\end{equation}
(cf.(\ref{e20})).

Now take an algebra generated by three generators
\begin{equation}
\label{e24}
ab-pba = F(N)\ ,\ qNa-aN = -a\ ,\  Nb-qbN = b\ ,
\end{equation}
(cf.(5)), where $p, q \in C,\ F$ a holomorphic function in a neighbourhood
of $\it Spec \{ N \} \cup \{ 0 \}$. 

{\bf \it THEOREM 6.}
For any $p,q \neq 0 \in {\bf C}$ in (23)  the identities (6.1), (6.2) and (10.1) -- (10.3) hold,
and also
\begin{equation}
\label{e25}
ba \left(ba - F \left({{N -1}\over q} \right) \right) \ldots 
\left( ba-\sum_{k=1}^n p^{k-1}F \left( {{N-\{ k \}}\over q^k}\right)\right) \ =
\ p^{n(n+1)\over2} b^{n+1} a^{n+1} \ .
\end{equation}
 
It is worth mentioning that for certain cases the formula (\ref{e25}) was known:
 $p=1$ \cite{4}, $p=q$ \cite{8}.

The proof is carried out by induction and is based on a certain generalization of Lemma 2 :

{\bf \it LEMMA 4.}
For any $p,q \neq 0 \in {\bf C}$ in (23)  and $n$ a natural number  
\renewcommand{\theequation}{25.{\arabic{equation}}}
\setcounter{equation}{0}
\begin{equation}
\label{e26.1}
a^n b^n \ =\ \sum_{k=0}^n \alpha_k (N) (ab)^k \ ,
\end{equation}
\begin{equation}
\label{e26.2}
b^n a^n \ =\ \sum_{k=0}^n \beta_k (N) (ab)^k \ ,
\end{equation}
where $\alpha, \beta$ are calculable functions.

\renewcommand{\theequation}{\arabic{equation}}
\setcounter{equation}{25}

We also use the simple observation that in (23)
\[
a f(N)= f(qN+1)a \quad , \quad  f(N)b= bf(qN+1) \ ,
\]
\[
[ f(N), a^nb^n] \ =\ [ f(N), b^na^n] \ =\ 0\ , \ n=1,2,\ldots \ ,
\]
where $f(N)$ is a holomorphic function in a neighbourhood
of $\it Spec \{ N \} \cup \{ 0 \}$.

It is worth noting that the algebra (5) or (13) can be interpreted as a 
deformation of the Heisenberg algebra. In turn, the algebra (\ref{e24})
contains as special cases: 

(i) $sl_2$ ( $p,q=1, F=2N$),

(ii) the $q$-deformed algebra $sl_{2q}$ ( $p=q^2, F=(q+1)N$, see (17)),

(iii) the quantum group $U_q(sl_2)$ ( $p,q=1, F= {\sin \tilde q
N \over { \tilde q- \tilde q ^{-1}}})$, \ and even,

(iv) the Heisenberg algebra  ( $q=1, F=1$) as a sub-algebra, \ and 

(v) the $q$-deformed Heisenberg algebra  ( $F=1$) as a sub-algebra. 

All the above demonstrates the general nature of the invariant identities (6.1), (6.2) and 
(10.1) -- (10.3).

\end{document}